\documentclass[12pt]{amsart}
\usepackage{latexsym}
\usepackage{amssymb,amscd,amsthm,amsfonts,verbatim,amsopn,color,graphics}
\usepackage[dvips]{epsfig} 
\def\Z{{\mathbb Z}}

\def\C{{\mathbb C}}

\def\R{{\mathbb R}}
\def\N{{\mathbb N}}
\def\K{{\mathbb K}}
\def\P{{\mathbb P}}

\def\AA{{\mathcal A}}

\def\DD{{\mathcal D}}
\def\FF{{\mathcal F}}
\def\GG{{\mathcal G}}
\def\HH{{\mathcal H}}
\def\II{{\mathcal I}}
\def\LL{{\mathcal L}}
\def\MM{{\mathcal M}}
\def\NN{{\mathcal N}}

\def\Ss{{\mathcal S}}

\def\mfS{{\mathfrak S}}
\def\mfL{{\mathfrak L}}

\newcommand{\at}{\mathfrak A}

\newtheorem{thm}{Theorem}[section]
\newtheorem{df} [thm]{Definition}

\newtheorem{prop}[thm]{Proposition}

\newtheorem{rem}[thm]{Remark}

\newtheorem{expl}[thm]{Example}
\newtheorem{cstr}[thm]{Construction}
\numberwithin{equation}{section}

\newenvironment{explrm}{\begin{expl} \rm}{\end{expl}}

\newenvironment{dfrm}{\begin{df} \rm}{\end{df}}
\newenvironment{cstrrm}{\begin{cstr} \rm}{\end{cstr}}

\begin{document}

\title[De~Concini-Procesi wonderful arrangement models]
{De~Concini-Procesi wonderful \\ arrangement models\\[0.3cm]
A discrete geometer's point of view}

\author{Eva Maria Feichtner}

\address{
Department of Mathematics, ETH Zurich, 8092 Zurich, Switzerland}
\email{feichtne@math.ethz.ch}

\begin{abstract}
  This expository article outlines the construction of
  De~Concini-Procesi arrangement models and describes recent progress
  in understanding their significance from the algebraic, geometric,
  and combinatorial point of view. Throughout the exposition, a strong
  emphasis is given to combinatorial and discrete geometric data that
  lies at the core of the construction.
\end{abstract}


\maketitle

\tableofcontents


\section{An invitation to arrangement models}
\label{sect_intr}

The complements of coordinate hyperplanes in a real or complex vector
space are easy to understand: The coordinate hyperplanes in $\R^n$
dissect the space into $2^n$ open orthants; removing the coordinate
hyperplanes from $\C^n$ leaves the complex torus~$(\C^*)^n$. Arbitrary
subspace arrangements, i.e., finite families of linear subspaces, have
complements with far more intricate combinatorics in the real case,
and far more intricate topology in the complex case.
Arrangement models improve this complicated situation locally --
constructing an arrangement model means to alter the ambient space so
as to preserve the complement and to replace the arrangement by a divisor 
with normal crossings,
i.e., a collection of smooth hypersurfaces which locally intersect
like coordinate hyperplanes.  Almost a decade ago, De~Concini and
Procesi have provided a canonical construction of arrangement models --
{\em wonderful\/} arrangement models -- that had significant impact in
various fields of mathematics.

\begin{quote}
{\em Why should a discrete geometer be interested in this model 
construction?\/}
\end{quote}

\smallskip
Because there is a wealth of {\em wonderful\/} combinatorial and
discrete geometric structure lying at the heart of the matter. Our aim
here is to bring these discrete pearls to light.

\medskip

First, combinatorial data plays a descriptive role at various places:
The combinatorics of the arrangement fully prescribes the model
construction and a natural stratification of the resulting space. We
will see details and examples in Section~\ref{sect_models}. In fact,
the rather coarse combinatorial data reflects enough of the situation
so as to, for instance, determine algebraic-topological invariants of
the arrangement models (compare the topological interpretation of the
algebra $D(\LL,\GG)$ that we study in Section~\ref{ssect_geomintDLG}).

Secondly, the combinatorial data that is put forward in the study of
arrangement models invites purely combinatorial generalizations. We
discuss such generalizations in Section~\ref{sect_comb} and show in the
subsequent Section~\ref{sect_btgeom} how this combinatorial generalization 
opens rather unexpected views when related back to geometry.

Finally, we propose arrangement models as a tool for resolving group
actions on manifolds in Section~\ref{sect_gractns}. 
Again, it is the open eye for discrete core data that enables the construction.

\medskip

We have attempted to keep the exposition rather self-contained and to
illustrate the development with many examples. We invite discrete
geometers to discover an algebro-geometric context in which familiar
discrete structures play a key role. We hope that yet many more
bridges will be built between algebraic and discrete geometry -- areas
that, despite the differences in terminology, concepts, and methods,
share what has inspired and driven mathematicians for centuries: the
passion for geometry.

\medskip

{\bf Acknowledgments:} I wish to thank the organizers of the workshop  
``Combinatorial and Discrete Geometry,'' held at MSRI in November 2003, 
for creating an event of out-most breadth, a true kaleidoscope
of topics unified by the ubiquity of geometry {\em and} combinatorics.


\section{Introducing the main character} 
\label{sect_models}

\noindent
We start out with explaining the De~Concini-Procesi arrangement model
construction. We will study some simple examples, which are rich
enough to convey the essential features of the models. Moreover, we
will outline some of the background and motivation for the model
construction.

\subsection{Basics on arrangements}

We first need to fix some basic terminology, in particular as it
concerns the combinatorial data of an arrangement. We suggest that the
reader, who is not familiar with the setting, reads through the first
part of this Section and compares the notions to the illustrations given 
for braid arrangements in Example~\ref{ex_braid}.

An {\em arrangement\/} $\AA\,=\,\{U_1,\ldots, U_n \}$ is a finite
family of linear subspaces in a real or complex vector space~$V$.  The
topological space associated first hand to such an arrangement is its
{\em complement\/} in the ambient space, 
$\MM(\AA):=V\,{\setminus}\,\bigcup\,\AA$.

Having arrangements in real vector space in mind, the topology of
$\MM(\AA)$ does not look very interesting: the complement is a
collection of open polyhedral cones, other than their number there is
no significant topological data connected to it. In the complex case,
however, already a single ``hyperplane'' in~$\C^1$,
i.e. $\AA\,{=}\,\{0\}$, has a nontrivial complement: it is homotopy
equivalent to $S^1$, the $1$-dimensional sphere.  The complement of
two (for instance, coordinate) hyperplanes in $\C^2$ is homotopy
equivalent to the torus $S^1{\times}\,S^1$.

The combinatorial data associated with an arrangement is customarily
recorded in a partially ordered set, the {\em intersection lattice\/}
$\LL\,{=}\,\LL(\AA)$. It is the set of intersections of subspaces in
$\AA$ ordered by reversed inclusion. We adopt terminology from the
theory of partially ordered sets and often denote the unique minimum
in $\LL(\AA)$ (corresponding to the empty intersection, i.e., the
ambient space $V$) by $\hat 0$ and the unique maximum of $\LL(\AA)$
(the overall intersection of subspaces in $\AA$) by $\hat 1$.  In many
situations, the elements of the intersection lattice are labeled by
the codimension of the corresponding intersection. For arrangements of
hyperplanes, this information is recorded in the rank function of the
lattice - the codimension of an intersection $X$ is the number of
elements in a maximal chain in the half-open interval $(\hat 0,X]$ 
in~$\LL(\AA)$.

As with any poset, we can consider the {\em order complex\/} 
$\Delta(\overline{\LL})$ of the proper part, 
$\overline{\LL}\,{:=}\,\LL\setminus\{\hat 0,\hat 1\}$, 
of the intersection lattice, i.e., the abstract simplicial complex 
formed by the linearly ordered subsets in $\overline{\LL}$,  
\[
     \Delta(\overline{\LL})\, = \, 
     \{X_1<\ldots <X_k\, | X_i \in \LL \setminus \{\hat 0,\hat 1\} \}\, .
\]
The topology of  $\Delta(\overline{\LL})$ plays a prominent role for 
describing the 
topology of arrangement complements. For instance, it is the crucial 
ingredient for the explicit description of cohomology groups of $\MM(\AA)$
by Goresky and MacPherson~\cite[Part III]{GM}. 

For hyperplane arrangements, the homotopy type of
$\Delta(\overline{\LL})$ is well-known: the complex is homotopy
equivalent to a wedge of spheres of dimension equal to the codimension
of the total intersection of $\AA$. The number of spheres can as well
be read from the intersection lattice, it is the absolute value of its
M\"obius function. For subspace arrangements however, the barycentric
subdivision of any finite simplicial complex can appear as the order
complex of the intersection lattice.

Besides $\Delta(\overline{\LL})$, we will often refer to the cone over
$\Delta(\overline{\LL})$ obtained by extending the linearly ordered
sets in $\overline{\LL}$ by the maximal element $\hat 1$ in $\LL$. We
will denote this complex by $\Delta(\LL{\setminus}\{\hat 0\})$ or
$\Delta(\LL_{>\hat 0})$.

\medskip
In order to have a standard example at hand, we briefly discuss braid
arrangements. This class of arrangements has figured prominently in
many places and has helped develop lots of arrangement theory
over the last decades.

\begin{explrm} 
\label{ex_braid}{\bf (Braid arrangements)} \newline
The arrangement $\AA_{n-1}$  given by 
the hyperplanes
\[
    H_{ij}: \,\, x_i=x_j\, , \quad \quad \mbox{for }\, 1\leq i<j\leq n\, ,
\]
in real $n$-dimensional vector space
is called the {\em (real) rank $n{-}1$ braid arrangement\/}. 
There is a complex version of this arrangement. It consists of hyperplanes 
$H_{ij}$ in $\C^n$ given by the {\em same\/} linear equations. We denote
the arrangement by $\AA_{n-1}^{\C}$. Occasionally, we will use the analogous
$\AA_{n-1}^{\R}$ if we want to stress the real setting.
In many situations a similar reasoning applies to the real and to the 
complex case. To simplify notation, we then use $\K$ to denote $\R$ or $\C$.

\begin{figure}[h]
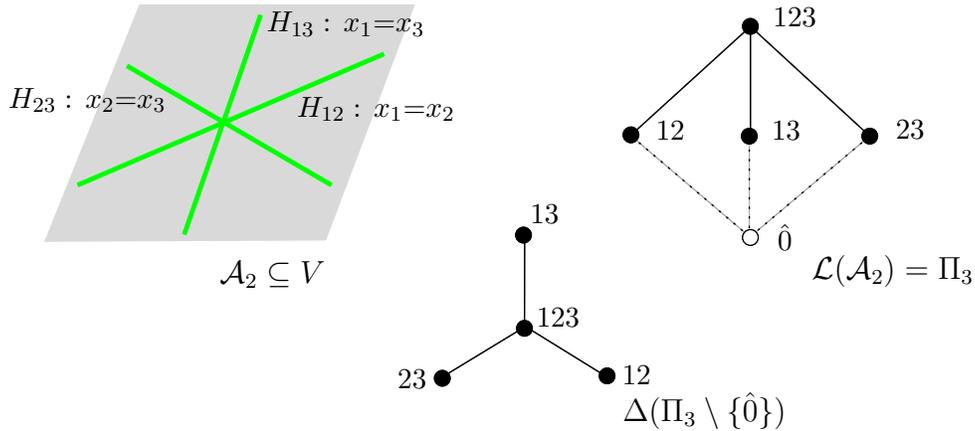

  \begin{picture}(0,0)%
    \includegraphics{A2.pstex}%
  \end{picture}%
  \input{A2.pstex_t}%
  
\caption{The rank $2$ braid arrangement $\AA_2$, 
its intersection lattice~$\Pi_3$,
and the order complex $\Delta(\Pi_3{\setminus}\,{\{\hat 0\}})$.}
\label{fig_A2}
\end{figure}

Observe that the diagonal $\Delta=\{x\in \K^n\,|\, x_1=\ldots =x_n\}$ 
is the overall intersection of hyperplanes in $\AA_{n-1}$. 
Without loosing any relevant information on the topology of the 
complement, we will often consider $\AA_{n-1}$ 
as an arrangement in complex or real $(n{-}1)$-dimensional space 
$V=\K^n/\Delta\,\cong\,$ $\{x\in \K^n\,|$ $\sum x_i=0\}$. 
This explains the,
at first sight, unusual indexing for braid arrangements.

The complement $\MM(\AA_{n-1}^{\R})$ is a collection of $n!$ polyhedral
cones, corresponding to the $n!$ linear orders on $n$ pairwise non-coinciding
coordinate entries.  The complement $\MM(\AA_{n-1}^{\C})$ is the classical
configuration space of the complex plane 
\[
   F(\C,n) \, \, = \, \, \{(x_1,\ldots,x_n)\in \C^n\,|\,
                   x_i\neq x_j \mbox{ for }\, i\neq j\,\}\, .
\]
This space is the classifying space of the pure braid group, which
explains the occurrence of the term ``braid'' for this class of arrangements.

As the intersection lattice of the braid arrangement $\AA_{n-1}$ we
recognize the {\em partition lattice\/}~$\Pi_n$, i.e., the set of set
partitions of $\{1,\ldots,n\}$ ordered by reversed refinement. The
correspondence to intersections in the braid arrangement can be easily
described: The blocks of a partition correspond to sets of coordinates
with identical entries, thus to the set of points in the corresponding
intersection of hyperplanes.

The order complex $\Delta(\overline{\Pi_n})$ is  a pure, $(n{-}1)$-dimensional 
complex that is homotopy equivalent to a wedge of $(n{-}1)!$ spheres of 
dimension~$n{-}1$.

\medskip
In Figure~\ref{fig_A2} we depict the real rank~$2$ braid arrangement
$\AA_2$ in $V\,{=}\,\R^3/\Delta$, its intersection lattice $\Pi_3$, and the
order complex $\Delta(\Pi_3\,{\setminus}\,\{\hat 0\})$. We denote
partitions in $\Pi_3$ by their non-trivial blocks. The depicted complex
is a cone over $\Delta(\overline{\Pi_3})$, a union of three points,
which indeed is the wedge of two $0$-dimensional spheres.
\end{explrm}


\subsection{The model construction}
\label{sect_DPmodels}

We provide two alternative definitions for De~Concini-Procesi
arrangement models: the first one describes the models as closures of
open embeddings of the arrangement complements. It comes in handy
for technical purposes.  Much more intuitive and suitable for
inductive constructions and proofs is the second definition, which 
describes arrangement models as results of sequences of blowups.

\begin{dfrm}
\label{df_model1}
{\bf (Model construction I)}\newline
Let $\AA$ be an arrangement of real or complex linear subspaces in~$V$.
Consider the map
\begin{eqnarray}
\Psi: \quad \MM(\AA) & \longrightarrow & 
                    V\, \times 
                    \, \prod_{X\in\LL_{>\hat 0}}\, \P(V/X) \label{eq_Psi} \\ 
       x & \longmapsto &   
           (\,x\,, (\langle x,X \rangle/X)_{X\in \LL_{> \hat 0}})\, ; 
\nonumber
\end{eqnarray}
it encodes the relative position of each point in the arrangement 
complement $\MM(\AA)$ with respect to the intersection of subspaces in $\AA$.
The map $\Psi$ is an 
open embedding; the closure of its image 
is called the {\em (maximal) De~Concini-Procesi wonderful model 
for $\AA$\/} and is denoted by~$Y_{\AA}$. 
\end{dfrm}

\begin{dfrm}
\label{df_model2}
{\bf (Model construction II)}\newline
Let $\AA$ be an arrangement of real or complex linear subspaces in~$V$.
Let $X_1,\ldots,X_t$ be a linear extension of the opposite order 
$\LL_{>\hat 0}^{\rm op}$ on~$\LL_{>\hat 0}$. 
The {\em (maximal) De~Concini-Procesi wonderful model 
for $\AA$\/}, $Y_{\AA}$, is the result of successively blowing up 
subspaces $X_1,\ldots,X_t$, respectively their proper transforms. 
\end{dfrm}

\medskip To avoid confusion with spherical blowups that have been
appearing in model constructions as well~\cite{G2}, let us emphasize
here that, also in the real setting, we think about blowups as
substituting points by projective spaces. Before we list the main
properties of arrangement models let us look at a first example.

\begin{explrm} 
\label{ex_YA2}
{\bf (The arrangement model $Y_{\AA_2}$)}\newline 
We consider the rank~$2$ braid arrangement $\AA_2$ in $V\,{=}\,\R^3/\Delta$.
Following the description in Definition~\ref{df_model2} we obtain 
$Y_{\AA_2}$ by
a single blowup of $V$ at $\{0\}$. The result is an open M\"obius band;
the exceptional divisor $D_{123}\,{\cong}\,\R\P^1$ in $Y_{\AA_2}$ intersects
transversally with the proper transforms $D_{ij}$ of the hyperplanes
$H_{ij}$, $1\,{\leq}\,i\,{<}\,j\,{\leq}\,3$. 
We illustrate the blowup in Figure~\ref{fig_YA2}.  

\begin{figure}[ht]
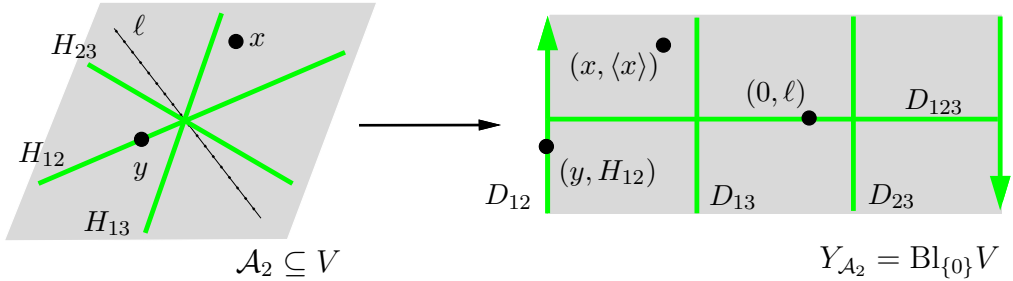

  \begin{picture}(0,0)%
    \includegraphics{YA2.pstex}%
  \end{picture}%
  \input{YA2.pstex_t}%
  
\caption{The maximal wonderful model for $\AA_2$.} \label{fig_YA2}
\end{figure}

In order to recognize the M\"obius band as the
closure of the image of $\Psi$ according to
Definition~\ref{df_model1}, observe that the product on the right-hand
side of~(\ref{eq_Psi}) consists of two relevant factors,
$V\,{\times}\,\R\P^1$. A point $x$ in $\MM(\AA_2)$ gets mapped to
$(x,\langle x \rangle)$ and we observe a one-to-one correspondence between
points in $\MM(\AA_2)$ and points in $Y_{\AA_2}\,{\setminus}\,
(D_{123}\,{\cup}\,D_{12}\,{\cup}\,D_{13} \,{\cup}\,D_{23})$. Points
which are added when taking the closure are of the form $(y,H_{ij})$
for $y\,{\in}\,H_{ij}{\setminus}\{0\}$ and $(0,\ell)$ for~$\ell$ some
line in~$V$.

Observe that the triple intersection of hyperplanes in $V$ has been
replaced by double intersections of hypersurfaces in $Y_{\AA_2}$.
Without changing the topology of the arrangement complement, the
arrangement of hyperplanes has been replaced by a normal crossing
divisor. Moreover, note that the irreducible divisor components
$D_{12}$, $D_{13}$, $D_{23}$, and $D_{123}$ intersect if and only if
their indexing lattice elements form a chain in $\LL(\AA_2)$.
\end{explrm}

The observations we made for $Y_{\AA_2}$ are special cases of the main
properties of (maximal) De~Concini-Procesi models that we list in the
following:

\begin{thm} \label{thm_models}
{\rm \cite[3.1 Thm., 3.2 Thm.]{DP1}}\newline
{\bf (1)} The arrangement model $Y_{\AA}$ as defined in
\ref{df_model1} and \ref{df_model2} is a smooth variety with a natural
projection map to the original ambient space, $\pi:\,
Y_{\AA}\longrightarrow V$, which is one-to-one on the arrangement
complement $\MM(\AA)$. \newline 
{\bf (2)} The complement of $\pi^{-1}(\MM(\AA))$ in $Y_{\AA}$ 
is a divisor with normal crossings; its irreducible components are the 
proper transforms $D_X$ of intersections $X$ in~$\LL$,
\[
      Y_{\AA} \setminus \pi^{-1}(\MM(\AA))\, \, =\, \, 
                        \bigcup_{X\in \LL_{>\hat 0}} D_X \, .
\]
{\bf (3)} Irreducible components $D_X$ for
$X\,{\in}\,\Ss\,{\subseteq}\,\LL_{>\hat 0}$ intersect if and only
if $\Ss$ is a linearly ordered subset in $\LL_{>\hat 0}$. If we
think about $Y_{\AA}$ as stratified by the irreducible components
of the normal crossing divisor and their intersections, then the poset
of strata coincides with the face poset of the order complex
$\Delta(\LL_{>\hat 0})$.
\end{thm}

\begin{explrm} 
\label{ex_A3}
{\bf (The arrangement model $Y_{\AA_3}$)}\newline Let us now consider
a somewhat larger and more complicated example, the rank~$3$ braid
arrangement $\AA_3$ in $V\,{\cong}\,\R^4/\Delta$. First note that the
intersection lattice of $\AA_3$ is the partition lattice $\Pi_4$,
which we depict in Figure~\ref{fig_Pi4} for later reference. Again, 
we denote partitions by their non-trivial blocks. 

\begin{figure}[ht]
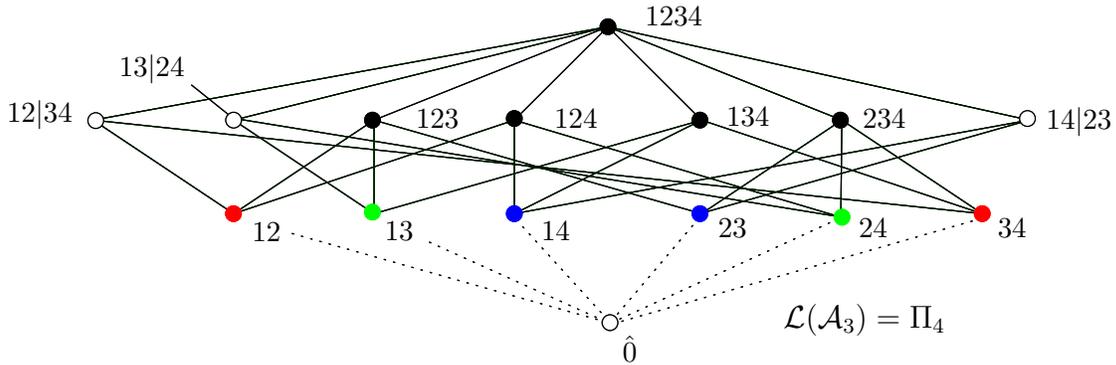

  \begin{picture}(0,0)%
    \includegraphics{Pi4.pstex}%
  \end{picture}%
  \input{Pi4.pstex_t}%
  
\caption{The intersection lattice of $\AA_3$.} \label{fig_Pi4}
\end{figure}

Following again the description of arrangement models given in
Definition~\ref{df_model2}, the first step is to blow up $V$ at
$\{0\}$. We obtain a line bundle over $\R\P^2$; in
Figure~\ref{fig_BlA3} we depict the exceptional divisor $D_{1234}
\cong \R\P^2$ stratified by the intersections of proper transforms of
hyperplanes in $\AA_3$.

\begin{figure}[ht]
  \begin{picture}(0,0)%
    \includegraphics{BlA3.pstex}%
  \end{picture}%
  \input{BlA3.pstex_t}%
  
\caption{The construction of $Y_{\AA_3}$.} \label{fig_BlA3}
\end{figure}

This first step is now followed by the blowup of triple, respectively 
double intersections of proper transforms of hyperplanes in arbitrary order. In each such
intersection the situation locally corresponds to the blowup of a
2-dimensional real vector space in a point as discussed in
Example~\ref{ex_YA2}. Topologically, the arrangement model $Y_{\AA_3}$
is a line bundle over a space obtained from a $7$-fold punctured
$\R\P^2$ by gluing $7$ M\"obius bands along their boundaries into the
boundary components.

We can easily check the statements of Theorem~\ref{thm_models} for 
$Y_{\AA_3}$. In particular, we see that intersections of 
irreducible divisors in $Y_{\AA_3}$ are non-empty if and only if the 
corresponding index sets form a chain in $\LL_{>\hat 0}$.
For instance, the $0$-dimensional stratum of the divisor stratification 
that is encircled in Figure~\ref{fig_BlA3} corresponds to the 
chain $14<134<1234$ in $\Pi_4\,{\setminus}\,\{\hat 0\}$. 
For comparison, we depict the order complex of 
$\Pi_4\,{\setminus}\,\{\hat 0\}$ in Figure~\ref{fig_NPi4max}. 
Recall that the complex 
is a pure $2$-dimensional cone with apex $1234$ over 
$\Delta(\overline{\Pi_4})$; we only draw its base.

\begin{figure}[ht]
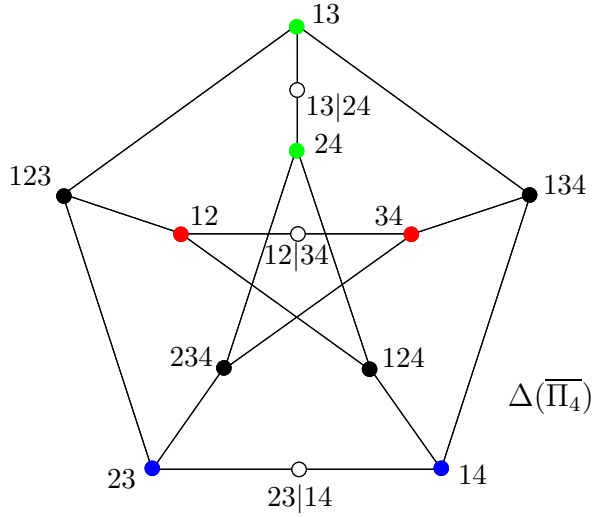

  \begin{picture}(0,0)%
    \includegraphics{NPi4max.pstex}%
  \end{picture}%
  \input{NPi4max.pstex_t}%
  
\caption{The order complex $\Delta(\overline{\Pi_4})$.} \label{fig_NPi4max}
\end{figure}

If our only objective was to construct a model for $\MM(\AA_3)$ with a
normal crossing divisor, it would be enough to blow up ${\rm
  Bl}_{\{0\}}V$ in the $4$ triple intersections. The result would be a
line bundle over a $4$-fold punctured~$\R\P^2$ with $4$ M\"obius bands
glued into boundary components.

This observation leads to a generalization of the model construction
presented so far: it is enough to do successive blowups on a specific
{\em subset\/} of intersections in $\AA$ to obtain a model with similar
properties as those summarized in Theorem~\ref{thm_models}.  In fact,
appropriate subsets of intersections lattices, so-called {\em building
sets\/}, were specified in~\cite{DP1}; all give rise to wonderful
arrangement models in the sense of Theorem~\ref{thm_models}.  The only
reservation being that the order complex $\Delta(\LL_{>\hat 0})$ is no
longer indexing non-empty intersections of irreducible divisors:
chains in $\LL_{>\hat 0}$ are replaced by so-called {\em nested
sets\/} -- subsets of building sets that again form an abstract
simplicial complex.

We will not give the original definitions of De~Concini and Procesi 
for building sets and nested sets in
this survey.  Instead, we will present a generalization of these notions 
for arbitrary meet-semilattices in Section~\ref{ssect_combres}. This
combinatorial abstraction has proved useful in many cases beyond
arrangement model constructions. Its relation to the original geometric 
context will be explained in Section~\ref{ssect_strat}.
\end{explrm}


\subsection{Some remarks on history} 
\label{ssect_hist}

Before we proceed, we briefly sketch the historic background of
De~Concini-Procesi arrangement models. Moreover, we outline an
application to a famous problem in arrangement theory that, among
other issues, served as a motivation for the model construction.

\medskip

Compactifications of configuration spaces due to 
Fulton and MacPherson~\cite{FuM} have prepared the scene for 
wonderful arrangement models. 
Their work is concerned with classical 
configurations spaces $F(X,n)$ of smooth algebraic varieties~$X$, i.e., 
spaces of $n$-tuples of pairwise distinct points in~$X$:
\[
     F(X,n) \,\,  = \,\, \{\, (x_1,\ldots, x_n) \in X^n\,| \, 
                  x_i\neq x_j \mbox{ for }\, i\neq j   \}\, . 
\]
A compactification~$X[n]$ of $F(X,n)$ is constructed in which the
complement of the original configuration space is a normal crossing
divisor; in fact, $X[n]$ has properties analogous to those listed for
arrangement models in~Theorem~\ref{thm_models}.  The relation to the
arrangement setting can be summarized by saying that, on the one hand,
the underlying spaces in the configuration space setting 
are incomparably more complicated -- smooth
algebraic varieties~$X$ rather than real or complex linear space; the
combinatorics, on the other hand, is far simpler -- it is the
combinatorics of our basic Examples~\ref{ex_YA2} and~\ref{ex_A3}, the
partition lattice~$\Pi_n$. The notion of building sets and nested sets, which
constitutes the defining combinatorics of arrangement models, has its 
roots in the Fulton-MacPherson construction for configuration spaces,
hence is inspired by the combinatorics of~$\Pi_n$.

Looking along the time line in the other direction, De~Concini-Procesi
arrangement models have triggered a number of more general
constructions with similar spirit: compactifications of conically
stratified complex manifolds by MacPherson and Procesi~\cite{MP}, and
model constructions for mixed real subspace and halfspace arrangements
and real stratified manifolds by Gaiffi~\cite{G2} that use spherical 
rather than classical blowups.

\medskip As a first impact, the De~Concini-Procesi model construction
has yielded substantial progress on a longstanding open question in
arrangement theory~\cite[Sect.\ 5]{DP1}, the question being whether
combinatorial data of a complex subspace arrangement determines the
cohomology algebra of its complement.  For arrangements of
hyperplanes, there is a beautiful description of the integral
cohomology algebra of the arrangement complement in terms of the
intersection lattice -- the {\em Orlik-Solomon algebra\/}~\cite{OS}.
Also, a prominent application of Goresky and MacPherson's Stratified
Morse Theory states that cohomology of complements of (complex and
real) subspace arrangements, as graded groups over $\Z$, are
determined by the intersection lattice and its codimension labelling.
In fact, there is an explicit description of cohomology groups in
terms of homology of intervals in the intersection lattice~\cite[Part
III]{GM}.  However, whether {\em multiplicative\/} structure is
determined as well remained an open question 20 years after it had
been answered for arrangements of hyperplanes~(see~\cite{FZ,dL} for
results on particular classes of arrangements).

The De~Concini-Procesi construction allows to apply Morgan's theory on
rational models for complements of normal crossing divisors~\cite{M}
to arrangement complements and to conclude that their {\em rational\/}
cohomology algebras indeed are determined by the combinatorics of the
arrangement.  A key step in the description of the Morgan model is the
presentation of cohomology of divisor components and their
intersections in purely combinatorial terms~\cite[5.1, 5.2]{DP1}. For
details on this approach to arrangement cohomology,
see~\cite[5.3]{DP1}.

Unfortunately, the Morgan model is fairly complicated even for small
arrangements, and the approach is bound to rational coefficients. The
model has been considerably simplified in work of
Yuzvinsky~\cite{Y2,Y3}. In \cite{Y2} explicit
presentations of cohomology algebras for certain classes of
arrangements were given. However, despite an explicit conjecture of 
an integral model for
arrangement cohomology in~\cite[Conj.6.7]{Y2}, extending the result
to integral coefficients remained out of reach.  Only years later,
the question has been fully settled to the positive in work of Deligne,
Goresky and MacPherson~\cite{DGM} with a sheaf-theoretic approach, and
parallely by de~Longueville and Schultz~\cite{dLS} using rather elementary
topological methods: Integral cohomology algebras of complex 
arrangement complements are indeed determined by combinatorial data.


\section{The combinatorial core data - a step beyond geometry}
\label{sect_comb}

We will now abandon geometry for a while and in this section fully 
concentrate on combinatorial and algebraic gadgets that are inspired 
by De~Concini-Procesi arrangement models. 

We first present a combinatorial analogue of De~Concini-Procesi
resolutions on purely order theoretic level 
following~\cite[Sect.\ 2\&3]{FK1}.
Based on the notion of building sets and nested sets for arbitrary 
lattices proposed therein, we define a family of commutative graded 
algebras for any given lattice.

The next Section then will be devoted to relate these objects
to geometry -- to the original context of De~Concini-Procesi arrangement 
models and, more interestingly so, to different seemingly unrelated 
contexts in geometry.

\subsection{Combinatorial resolutions}
\label{ssect_combres}

We will state purely combinatorial definitions of {\em building sets\/} 
and {\em nested sets\/}. Recall that, in the context of model
constructions, building sets list the strata that are to be blown up in
the construction process, and nested sets describe beforehand the
non-empty intersections of irreducible
divisor components in the final resolution.

\medskip 

Let $\LL$ be a finite meet-semilattice, i.e., a finite poset such that
any pair of elements has a unique maximal lower bound. In particular, 
such a meet-semilattice has a unique minimal element that we denote with 
$\hat 0$. We will talk about semilattices for short. As a basic reference
on partially ordered sets we refer to~\cite[Ch.\ 3]{St}.

\begin{dfrm} \label{df_building}
{\bf (Combinatorial building sets)}\newline
A subset $\GG\,{\subseteq}\,\LL_{>\hat 0}$ 
in a finite meet-semilattice~$\LL$  is called a 
{\em building set\/} if for any 
$X\,{\in}\,\LL_{>\hat 0}$  and
{\rm max}$\, \GG_{\leq X}=\{G_1,\ldots,G_k\}$ there is an isomorphism
of posets
\begin{equation}\label{eq_buildg}
\varphi_X:\,\,\, \prod_{j=1}^k\,\,\, [\hat 0,G_j] \,\, 
                               \stackrel{\cong}{\longrightarrow}
                                \,\, [\hat 0,X]
\end{equation}
with $\varphi_X(\hat 0, \ldots, G_j, \ldots, \hat 0)\, = \, G_j$ 
for $j=1,\ldots, k$. 
We call $F_{\GG}(X)\,{:=}\, 
\mathrm{max}\,\GG_{\leq X}$  the {\em set of factors\/} of $X$ in $\GG$.
\end{dfrm}

There are two extreme examples of building sets for any semilattice:
we can take the full semilattice $\LL_{>\hat 0}$ as a building set. 
On the other hand, the set of elements $X$ in $\LL_{>\hat 0}$ which 
do not allow for a 
product decomposition of the lower interval~$[\hat 0,X]$ form the unique
minimal building set (see Example~\ref{ex_bsns} below). 

Intuitively speaking, building sets are formed by elements in the 
semilattice that are the perspective factors of product decompositions. 

\medskip
Any choice of a building set $\GG$ in $\LL$ gives rise to a family of
so-called {\em nested sets\/}. These are, roughly speaking, subsets of
$\GG$ whose antichains are sets of factors with respect to the chosen
building set.  Nested sets form an abstract simplicial complex
on the vertex set $\GG$. This simplicial complex plays
the role of the order complex for arrangement models more general than
the maximal models discussed in Section~\ref{sect_DPmodels}.

\begin{dfrm} 
\label{df_nested}
{\bf (Nested sets)} \newline
Let $\LL$ be a finite meet-semilattice and $\GG$ a building set in $\LL$.
A subset $\Ss$ in $\GG$ is called~{\em nested\/} (or $\GG$-{\em nested\/}
if specification is needed) 
if, for any set of incomparable elements 
$X_1,\dots,X_t$ in $\Ss$ of cardinality at least two, 
the join $X_1\vee\dots\vee X_t$ exists and does not belong to $\GG$.
The $\GG$-nested sets form an abstract simplicial complex $\NN(\LL,\GG)$,
the {\em nested set complex\/} with respect to $\LL$ and $\GG$. 
\end{dfrm}

Observe that if we choose the full semilattice as a building set, then
a subset is nested if and only if it is linearly ordered in
$\LL$. Hence, the nested set complex $\NN(\LL,\LL_{>\hat 0})$
coincides with the order complex $\Delta(\LL_{>\hat 0})$.

\begin{explrm} 
\label{ex_bsns}
{\bf (Building sets and nested sets for the partition lattice)} \newline
Choosing the maximal building set in the partition lattice $\Pi_n$, 
we obtain the order complex $\Delta((\Pi_n)\,{\setminus}\,\{\hat 0\})$ 
as the associated complex of nested sets. Topologically, it is 
a cone over a wedge of $(n-1)!$ spheres of dimension $n{-}1$. 

The minimal building set $\GG_{{\rm min}}$ in $\Pi_n$ is given by
partitions with at most one block of size larger or equal~$2$, the so-called 
modular elements in~$\Pi_n$.  We can
identify these partitions with subsets of $\{1,\ldots,n\}$ of size
larger or equal~$2$.  A collection of such subsets is nested, if and
only if none of the pairs of subsets have a non-trivial intersection,
i.e., for any pair of subsets they are either disjoint or one is
contained in the other. Referring to a naive picture of such
containment relation explains the choice of the term {\em nested} --
it appeared first in the work of Fulton and MacPherson~\cite{FuM} on
compactifications of classical configuration spaces.  As we noted earlier, 
the combinatorics they are concerned with is indeed the combinatorics of the
partition lattice.

For the rank~$3$ partition lattice $\Pi_3$, 
maximal and minimal building sets coincide, 
$\GG=\Pi_3\,{\setminus}\,\{\hat 0\}$. 
The nested set complex $\NN(\Pi_3,\GG)$
is the order complex $\Delta(\Pi_3\,{\setminus}\,\{\hat 0\})$ 
depicted in Figure~\ref{fig_A2}.

For the rank~$4$ partition lattice $\Pi_4$, we have seen the nested
set complex for the maximal building set $\NN(\Pi_4,\GG_{{\rm max}})$
in Figure~\ref{fig_NPi4max}.  The nested set complex associated with
the minimal building set $\GG_{{\rm min}}$ in $\Pi_4$ is depicted in
Figure~\ref{fig_NPi4min}.  Again, $\NN(\Pi_4,\GG_{{\rm min}})$ is a
cone with apex $1234$, and we only draw its base,
$\NN(\overline{\Pi_4},\GG_{{\rm min}})$.

\begin{figure}[ht]
  \begin{picture}(0,0)%
    \includegraphics{NPi4min.pstex}%
  \end{picture}%
  \input{NPi4min.pstex_t}%
  
\caption{The nested set complex $\NN(\overline{\Pi_4},\GG_{{\rm min}})$.} 
\label{fig_NPi4min}
\end{figure}

Adding one or two $2$-block partitions to $\GG_{{\rm min}}$ yields all
the other building sets for $\Pi_4$. The corresponding nested set
complexes are subdivisions of $\NN(\Pi_4,\GG_{{\rm min}})$.

When studying the (maximal) wonderful model $Y_{\AA_3}$ in
Example~\ref{ex_A3} we had observed that, if we only wanted to
achieve a model with normal crossing divisors, it would have been
enough to blow up the overall and the triple intersections.  This
selection of strata, respectively elements in
$\LL(\AA_3)\,{=}\,\Pi_4$, exactly corresponds to the minimal building
set $\GG_{{\rm min}}$ in~$\Pi_4$ -- a geometric motivation for
Definition~\ref{df_building}.

Let us also get a glimpse on the geometry that motivates
the definition of nested sets: comparing simplices in
$\NN(\Pi_4,\GG_{{\rm min}})$ with intersections of irreducible divisor
components in the arrangement model resulting from blowups along
subspaces in $\GG_{{\rm min}}$, we see that there is a $1$--$1$
correspondence. For instance, $\{12,34\}$ is a nested set with respect
to $\GG_{{\rm min}}$, and divisor components $D_{12}$ and $D_{34}$
intersect in the model (compare Figure~\ref{fig_BlA3}).
\end{explrm}

It is not a coincidence that, in the example above, one nested set
complex is a subdivision of the other if one building set contains the
other.  In fact, the following holds:

\begin{thm} 
\label{thm_nsc}
{\rm \cite[Prop.\ 3.3, Thm.\ 4.2]{FM}}
For any finite meet-semilattice $\LL$, and $\GG$ a building set in $\LL$, 
the nested set complex $\NN(\LL,\GG)$ is homotopy equivalent to the 
order complex of $\LL_{>\hat 0}$,
\[
\NN(\LL,\GG)\, \, \simeq \, \, \Delta(\LL_{>\hat 0})\, .
\]
Moreover, 
if $\LL$ is atomic, i.e., any element is a join of a set of atoms,
and $\GG$ and $\HH$ are building sets with $\GG\,{\supseteq}\,\HH$, then
the nested set complex $\NN(\LL,\GG)$ is obtained from $\NN(\LL,\HH)$ by
a sequence of stellar subdivisions. In particular, the complexes are 
homeomorphic.
\end{thm}

We now propose a construction on semilattices producing new semilattices: 
the {\em combinatorial blowup\/} of a semilattice in an element. 

\begin{dfrm}
\label{df_cblowup}
{\bf (Combinatorial blowup)} \newline
For a semilattice $\LL$ and an element $X$ in $\LL_{>\hat 0}$ we define 
a poset $(\mathrm{Bl}_X\LL,\prec)$ on the set of elements
\[
  \mathrm{Bl}_X\LL\,  = \,   
   \{\,Y\,|\, Y\in \LL, Y\not \geq X\} \, \cup  \,  
   \{\,Y^{\prime}\,|\, Y\in \LL, Y\not \geq X, \,\mathrm{and}\,
                   Y\vee X \,\, \mathrm{exists \,\, in}\,\LL \,\} \, .
\]
The order relation $<$ in $\LL$ determines the order relation 
$\prec$ within the two parts of~$\mathrm{Bl}_X\LL$ described above,
\[
\begin{array}{lll}
    Y \prec Z\, , & \mathrm{for} &  Y<Z \,\,\,\mathrm{in}\,\, \LL\, , \\
    Y^{\prime} \prec Z^{\prime}\, , & \mathrm{for} &  Y<Z \,\,\,\mathrm{in}\,\, \LL\, , \\
    \end{array}
\]
and additional order relations between elements of these two parts are 
defined by 
\[
\begin{array}{lll} 
  Y \prec Z^{\prime}\, , & \mathrm{for} &  Y<Z \,\,\,\mathrm{in}\,\, \LL\, , 
\end{array}
\]
where in all three cases it is assumed that $Y,Z \not \geq X$ in $\LL$.
We call $\mathrm{Bl}_X\LL$ the {\em combinatorial blowup\/} of 
$\LL$ in $X$. 
\end{dfrm}

\begin{figure}[ht]
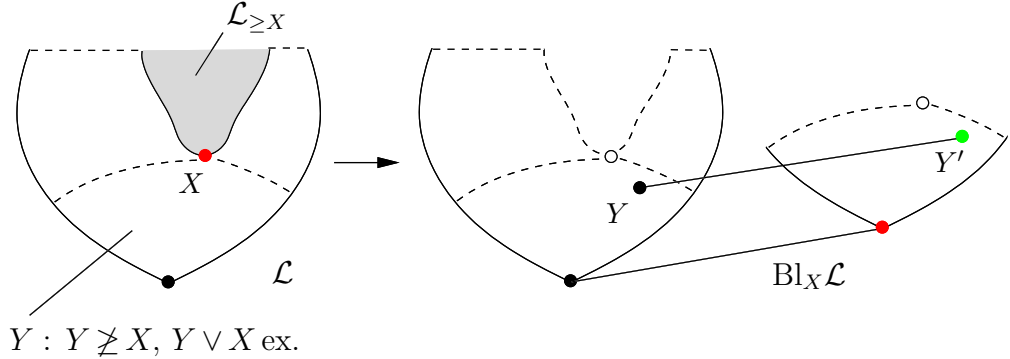

  \begin{picture}(0,0)%
    \includegraphics{combblowup.pstex}%
  \end{picture}%
  \input{combblowup.pstex_t}%
  
\caption{A combinatorial blowup.} 
\label{fig_combblowup}
\end{figure}

In fact, the poset $\mathrm{Bl}_X\LL$ is again a semilattice.
We believe that Figure~\ref{fig_combblowup} will much better explain what is going on.

The construction does the following: it removes the closed upper interval on
top of $X$ from $\LL$, and then marks the set of elements in $\LL$
that are not larger or equal $X$, but have a join with $X$ in
$\LL$. This subset of $\LL$ (in fact, a lower ideal in the sense of
order theory) is doubled and any new element $Y^{\prime}$ in the copy is
defined to be covering the original element $Y$ in $\LL$.  The order
relations in the remaining, respectively the doubled, part of $\LL$ 
stay the same as before.

\begin{figure}[ht]
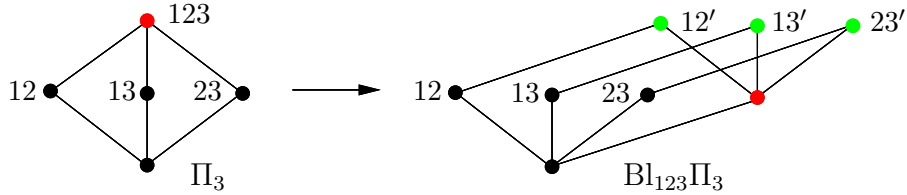

  \begin{picture}(0,0)%
    \includegraphics{BlPi3.pstex}%
  \end{picture}%
  \input{BlPi3.pstex_t}%
  
\caption{The combinatorial blowup of $\Pi_3$ in $123$.} 
\label{fig_BlPi3}
\end{figure}

In Figure~\ref{fig_BlPi3} we give a concrete example: the
combinatorial blowup of the maximal element $123$ in $\Pi_3$, ${\rm
  Bl}_{123} \Pi_3$.  The result should be compared with
Figure~\ref{fig_YA2}. In fact, ${\rm Bl}_{123} \Pi_3$ is the face
poset of the divisor stratification in $Y_{\AA_2}\,{=}\,{\rm
  Bl}_{\{0\}}V$.

The following theorem shows that the three concepts introduced above --
combinatorial building sets, nested sets, and combinatorial blowups --
fit together so as to provide a combinatorial analogue of the
De~Concini-Procesi model construction.

\begin{thm} \label{thm_combres}
{\rm \cite[Thm.\ 3.4]{FK1}}
Let $\LL$ be a semilattice, $\GG$ a combinatorial building set in
$\LL$, and $G_1,\ldots, G_t$ a linear order on $\GG$ that is
non-increasing with respect to the partial order on $\LL$. Then, 
consecutive combinatorial blowups in $G_1,\ldots, G_t$
result in the face poset of the nested set complex $\NN(\LL,\GG)$:
\[
   {\rm Bl}_{G_t}(\ldots({\rm Bl}_{G_2}({\rm Bl}_{G_1}\LL))\ldots)
   \, \, = \, \, \FF(\NN(\LL,\GG))\, . 
\]
\end{thm}

\subsection{An algebra defined for atomic lattices}
\label{ssect_algdef}

For any atomic lattice, we define a family of graded commutative algebras 
based on the notions of building sets and nested sets given above.
Our exposition here and in Section~\ref{ssect_geomintDLG} follows~\cite{FY}.
Restricting our attention to atomic lattices is not essential for the 
definition. For various algebraic considerations and for the geometric 
interpretations (cf.\ \ref{ssect_geomintDLG}), however, it is convenient 
to assume that the lattice is atomic. 

\begin{dfrm} 
\label{def_DLG}
  Let $\LL$ be a finite atomic lattice, $\at(\LL)$ its set of atoms, and
  $\GG$ a building set in~$\LL$. We define the 
  algebra $D(\LL,\GG)$ of $\LL$ with respect to $\GG$ as
\[ 
D(\LL,\GG) \, \, := \, \, 
             \Z\,[\{x_{G}\}_{G\in \GG}] \, \Big/ \, \II \, ,
\]
where the ideal of relations $\II$ is generated by 
\begin{eqnarray*}
    \prod_{i=1}^t \,x_{G_i} & & \qquad 
          \mbox{for }\,\,\{G_1,\ldots ,G_t\}\,\not\in \,\NN(\LL,\GG) \,,\\
    \sum_{G\geq H}\, x_G & & \qquad  
                            \mbox{for }\, H \in \at(\LL)\, .  
\end{eqnarray*}
\end{dfrm} 

Observe that this algebra is a quotient of the face ring 
of the nested set complex~$\NN(\LL,\GG)$. 

\begin{explrm} \label{ex_DPn}
{\bf (Algebras associated to $\Pi_3$ and $\Pi_4$) } \newline
For $\Pi_3$ and its only building set 
$\GG_{\rm max}\,{=}\, \Pi_3\,{\setminus}\,\{\hat 0\}$, the algebra 
reads as follows:

\begin{eqnarray*}
   D(\Pi_3,\GG_{\rm max}) &  = &  
       \Z\,[x_{12}, x_{13}, x_{23}, x_{123}] \Big/  
       \left\langle
         \begin{array}{l}
         x_{12}x_{13},\,\, x_{12}x_{23},\,\,  
         x_{13}x_{23}                               \\
         x_{12}+x_{123},\,\,  x_{13}+x_{123}, \,\, x_{23}+x_{123}
       \end{array}
       \right\rangle \\
& \cong & 
        \Z\,[x_{123}] / \langle x_{123}^2 \rangle\, .
\end{eqnarray*}
For  $\Pi_4$ and its minimal building set $\GG_{\rm min}$, we obtain 
the following algebra after slightly simplifying the presentation:
\begin{eqnarray*}
\lefteqn{ 
   D(\Pi_4,\GG_{\rm min}) \, \cong \, 
       \Z\,[x_{123}, x_{124}, x_{134}, x_{234}, x_{1234}] \, \Big/  }\\
              & &  \qquad \qquad \qquad \qquad
       \left\langle \begin{array}{ll} 
        x_{ijk} \,x_{1234} & \mbox{for all }\, 1{\leq}i{<}j{<}k{\leq}4 \\
        x_{ijk} \, x_{i'j'k'} & \mbox{for all }\, ijk \neq i'j'k' \\
        x_{ijk}^2 + x_{1234}^2 \quad & \mbox{for all }\, 
                                               1{\leq}i{<}j{<}k{\leq}4  
       \end{array}
       \right\rangle \, .
\end{eqnarray*}
\end{explrm}

There is an explicit description for a Gr\"obner basis of the ideal $\II$,
which in particular yields an explicit description for a monomial basis 
of the graded algebra $D(\LL,\GG)$.

\begin{thm}
\label{thm_Groebner}
{\bf (1)} {\rm \cite[Thm.\ 2]{FY}}
The following polynomials form a Gr\"obner basis of the ideal~$\II$:
\begin{eqnarray*}
& &  \prod_{G\in \Ss} \,x_{G} \qquad
                     \mbox{for }\,\, \Ss\not\in \NN(\LL,\GG)\, , \\
& &  \,\prod_{i=1}^k \,x_{A_i} \,  
                          \Big(\sum_{G\geq B}\, x_G \,\Big)^{d(A,B)}\, , 
\end{eqnarray*}
where $A_1,\ldots, A_k$ are maximal elements in a nested set $\HH\,{\in}\,
\NN(\LL,\GG)$, $B\,{\in}\,\GG$ with $B\,{>}\,A=\bigvee_{i=1}^k A_i$, and
$d(A,B)$ is the minimal number of atoms needed to generate $B$ from $A$ by 
taking joins. \newline 
{\bf (2)} {\rm \cite[Cor.\ 1]{FY}}
The resulting linear basis for the graded algebra $D(\LL,\GG)$
is given by the following set of monomials:
\[
    \prod_{A\in\Ss} x_A^{m(A)}\, ,
\]
where $\Ss$ is running over all nested subsets of $\GG$,
$m(A)\,{<}\,d(A',A)$, and $A'$ is the join of $\Ss\,{\cap}\,\LL_{<A}$.
\end{thm}

Part (2) of Theorem~\ref{thm_Groebner} 
generalizes a basis description by Yuzvinsky~\cite{Y1}  
for $D(\LL,\GG)$ in the case of $\GG$ being the minimal building set in an
intersection lattice $\LL$ of a complex hyperplane arrangement.  Yuzvinsky's
basis description has also been generalized in a somewhat different direction 
by Gaiffi~\cite{G1}, namely for closely related algebras associated with 
complex subspace arrangements. 

We will return to the algebra $D(\LL,\GG)$ and discuss its geometric 
significance in Section~\ref{ssect_geomintDLG}. 


\section{Returning to geometry}
\label{sect_btgeom}

\subsection{Understanding stratifications in wonderful models} 
\label{ssect_strat}

Let us first relate the combinatorial setting  of building sets and 
nested sets developed in Section~\ref{ssect_combres} to its origin, 
the De~Concini-Procesi 
model construction. Here is how to recover the original notion
of building sets~\cite[2.3 Def.]{DP1}, we call them {\em geometric building sets\/}, 
from our definitions:

\begin{dfrm} 
\label{df_geombuilding}
{\bf (Geometric building sets)} \newline
Let $\LL$ be the intersection lattice of an arrangement of subspaces 
in real or complex vector space $V$ and ${\rm cd}:\,\LL \rightarrow \N$ 
a function on $\LL$ assigning the codimension of the corresponding subspace 
to each lattice element. A subset $\GG$ in $\LL$ is a {\em geometric 
building set\/} if it is a building set in the sense of~\ref{df_building}, 
and for any $X\,{\in}\,\LL$ the codimension of~$X$ is equal to the sum of 
codimensions of its factors, $F_{\GG}(X)$:
\[
          {\rm cd}\, (X) \, = \, \sum_{Y\in F_{\GG}(X)} {\rm cd}\, (Y)\, .
\]
\end{dfrm}

An easy example shows that the notion of geometric building sets indeed
is more restrictive than the notion of combinatorial building sets. 
For arrangements of hyperplanes, however, the notions 
coincide~\cite[Prop.\ 4.5.(2)]{FK1}.

\begin{explrm} \label{ex_geombuilding}
{\bf (Geometric versus combinatorial building sets)} \newline
Let $\AA$ denote the following arrangement of $3$ subspaces in $\R^4$:
\[
A_1:\,\, x_4  =  0\, , \quad
A_2:\,\, x_1= x_2 =  0\, ,\quad
A_3:\,\, x_1= x_3 =  0 \,.
\]
The intersection lattice $\LL(\AA)$ is a boolean algebra on $3$
elements; we depict the lattice with its codimension labelling in
Figure~\ref{fig_geombuild}. The set of atoms obviously is a
combinatorial building set. However, any geometric building set must
contain the intersection $A_2\,{\cap}\, A_3$: its codimension is {\em
  not\/} the sum of codimensions of its (combinatorial) factors $A_2$
and $A_3$.

\begin{figure}[ht]
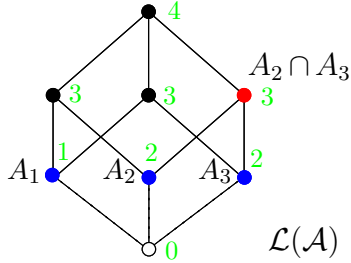

  \begin{picture}(0,0)%
    \includegraphics{geombuild.pstex}%
  \end{picture}%
  \input{geombuild.pstex_t}%
  
\caption{Geometric versus combinatorial building sets.} 
\label{fig_geombuild}
\end{figure}
\end{explrm}

\medskip 
As we mentioned before, there are wonderful model
constructions for arrangement complements $\MM(\AA)$ that start from
an arbitrary geometric building set $\GG$ of the intersection lattice
$\LL(\AA)$~\cite[3.1]{DP1}: In Definition~\ref{df_model1}, replace the
product on the right hand side of~(\ref{eq_Psi}) by a product over
building set elements in $\LL$, and obtain the wonderful model
$Y_{\AA,\GG}$ by again taking the closure of the image of $\MM(\AA)$
under $\Psi$. In Definition~\ref{df_model2}, replace the linear
extension of $\LL_{>\hat 0}^{\rm op}$ by a non-increasing linear order
on the elements in $\GG$, and obtain the wonderful model $Y_{\AA,\GG}$
by successive blowups of subspaces in $\GG$, and of proper transforms
of such.

The key properties of these models are analogous to those listed in
Theorem~\ref{thm_models}, where in part (2), lattice elements are
replaced by building set elements, and in part $(3)$, chains in $\LL$
as indexing sets of non-empty intersections of irreducible components
of divisors are replaced by nested sets.
Hence, the face poset of the stratification of
$Y_{\AA,\GG}$ given by irreducible components of divisors and their
intersections coincides with the face poset of the nested set complex
$\NN(\LL,\GG)$. Compare Examples~\ref{ex_A3} and \ref{ex_bsns}, where
we found that nested sets with respect to the minimal building set
$\GG_{{\rm min}}$ in $\Pi_4$ index non-empty intersections of
irreducible divisor components in the arrangement model 
$Y_{\AA_3,\GG_{{\rm min}}}$.

\medskip
While the intersection lattice $\LL(\AA)$ captures the combinatorics of 
the stratification of $V$ given by subspaces of $\AA$ and their 
intersections, the nested set complex $\NN(\LL,\GG)$ captures the 
combinatorics of the divisor stratification of the wonderful model 
$Y_{\AA,\GG}$. More than that: combinatorial blowups turn out to be the 
right concept to describe the incidence change of strata during the
construction of wonderful arrangement models by successive
blowups:

\begin{thm} {\rm \cite[Prop.\ 4.7 (1)]{FK1}}
Let $\AA$ be a complex subspace arrangement, $\GG$ a geometric building set 
in $\LL(\AA)$, and $G_1,\ldots, G_t$ a non-increasing linear order on $\GG$.
Let ${\rm Bl}_i(\AA)$ denote the result of blowing up strata 
$G_1,\ldots,G_i$, for 
$i\,{\leq}\,t$, and denote by $\LL_i$ the face poset of the stratification
of ${\rm Bl}_i(\AA)$ by proper transforms of subspaces in $\AA$ and the 
exceptional divisors. Then the poset  $\LL_i$ coincides with the successive  
combinatorial blowups of $\LL$ in $G_1,\ldots G_i$:
\[
   \LL_i \, \, = \, \,   
{\rm Bl}_{G_i}(\ldots({\rm Bl}_{G_2}({\rm Bl}_{G_1}\LL))\ldots) \, .           
\]  
\end{thm}

\medskip
Combinatorial building sets, nested sets and combinatorial blowups
occur in other situations and prove to be the right concept for
describing stratifications in more general model constructions. This
applies to the {\em wonderful conical compactifications\/} of
MacPherson and Procesi~\cite{MP} as well as to models for mixed
subspace and halfspace arrangements and for stratified real manifolds
by Gaiffi~\cite{G2}.

Also, combinatorial blowups describe the effect which stellar 
subdivisions in polyhedral fans have on the face  poset of the fans.  
In fact, combinatorial blowups describe the incidence change of torus orbits
for resolutions of toric varieties by consecutive blowups in closed
torus orbits. This implies, in particular, that for any toric variety
and for any choice of a combinatorial building set in the face poset
of its defining fan, we obtain a resolution of the variety with torus
orbit structure prescribed by the nested set complex associated to the
chosen building set. We believe that such combinatorially prescribed
resolutions can prove useful in various concrete 
situations~(see~\cite[Sect.4.2]{FK1} for further details).
 
\medskip


There is one more issue about nested set stratifications of maximal
wonderful arrangement models that we want to discuss here, mostly in
perspective of applications in Section~\ref{sect_gractns}. According
to Definition~\ref{df_model1}, any point in the model $Y_{\AA}$ can be
written as a collection of a point in $V$ and lines in $V$, one line for
each element in~$\LL(\AA)$. There is a lot of redundant information in
this rendering, e.g., points on the open stratum $\pi^{-1}(\MM(\AA))$ are fully
determined by their first ``coordinate entry,'' the point in 
$\MM(\AA)\,{\subseteq}\,V$. 

Here is a more economic encoding of a point  $\omega$ on 
$Y_{\AA}$~\cite[Sect.\ 4.1]{FK2}:
we find that $\omega$
can be uniquely written as
\begin{equation}\label{eq_encoding}
\omega \,\, = \,\, (x,H_1, \ell_1, H_2, \ell_2, \ldots, H_t,\ell_t)
        \,\, = \,\,(x,\ell_1, \ell_2, \ldots,\ell_t)\, , 
\end{equation}
where $x$ is a point in~$V$, the $H_1,\ldots,H_t$ form a descending chain 
of  subspaces in $\LL_{>\hat 0}$, and the~$\ell_i$ are lines in $V$.
More specifically, $x\,{=}\,\pi(\omega)$, and $H_1$ is the 
maximal lattice element that, as a subspace of~$V$,
contains~$x$. The line $\ell_1$ is orthogonal to $H_1$ and
corresponds to the coordinate entry of $\omega$ indexed by $H_1$
in $\P(V/H_1)$.  The lattice element $H_2$, in turn, is the
maximal lattice element that contains both $H_1$ and $\ell_1$. The
specification of lines~$\ell_i$, i.e., lines that correspond to
coordinates of $\omega$ in $\P(V/H_i)$, and the construction of
lattice elements $H_{i+1}$, continues analogously for $i\geq 2$
until a last line $\ell_t$ is reached whose span with $H_t$ is not
contained in any lattice element other than the full ambient
space~$V$. 
Observe that the $H_i$ are determined by $x$
and the sequence of lines $\ell_i$; we choose to include the $H_i$ 
in order to keep the notation more transparent.

The full coordinate information on $\omega$ can be recovered from 
(\ref{eq_encoding}) by setting $H_0=\bigcap \AA$, 
$\ell_0=\langle x\rangle$, and retrieving the coordinate $\omega_{H}$ 
indexed by $H\,{\in}\,\LL_{> \hat 0}$ as 
\[ 
\omega_H\, \, = \, \, \langle \ell_j,H\rangle / H  \, \, \in \, \, \P(V/H)\, ,
\]
where $j$ is chosen from $\{1,\ldots, t\}$ such that
$H\leq H_j$, but $H\not \leq H_{j+1}$.

A nice feature of this encoding is that for a given point $\omega$ in 
$Y_{\AA}$ we can tell the open stratum in the nested set 
stratification which contains it:

\begin{prop}
\label{prop_charopenstrat}
{\rm (\cite[Prop 4.5]{FK2})}
A point $\omega$ in a maximal arrangement model $Y_{\AA}$ is contained 
in the open stratum indexed by the chain
$H_1\,{>}\,H_2\,{>}\,\ldots \,{>}\,H_t$ in~$\LL_{>\hat 0}$ if 
and only if its point/line description~{\rm (\ref{eq_encoding})} reads 
$\omega\,{=}\,(x,H_1, \ell_1, H_2, \ell_2, \ldots, H_t,\ell_t)$.
\end{prop}


\subsection{A wealth of geometric meaning for $D(\LL,\GG)$}         
\label{ssect_geomintDLG}

We turn to the algebra $D(\LL,\GG)$ that we defined for any atomic
lattice $\LL$ and combinatorial building set $\GG$ in $\LL$ in
Section~\ref{ssect_algdef}. We give two geometric interpretations for
this algebra; one is restricted to $\LL$ being the intersection
lattice of a complex hyperplane arrangement and originally motivated
the definition of~$D(\LL,\GG)$, the other applies to any atomic
lattice and provides for a somewhat unexpected connection to toric
varieties.

Let us briefly comment on the projective version of wonderful
arrangement models that we need in this context (see~\cite[\S 4]{DP1}
for details).  For any arrangement of linear subspaces $\AA$ in $V$, a
model for its projectivization $\P \AA\,{=}\,\{\P A\,|\,
A\,{\in}\,\AA\}$ in $\P V$, i.e., for $\MM(\P \AA)\,{=}\,\P
V\,{\setminus}\, \bigcup \P \AA$, can be obtained by replacing the
ambient space $V$ by its projectivization $\P V$ in the model
constructions~\ref{df_model1} and~\ref{df_model2}. The constructions
result in a smooth projective variety that we denote
by~$Y_{\AA}^{\P}$.  A model $Y_{\AA,\GG}^{\P}$ for a specific
geometric building set $\GG$ in $\LL$ can be obtained analogously. In
fact, under the assumption that $\P(\bigcap \AA)$ is contained in the
building set~$\GG$, the affine model $Y_{\AA,\GG}$ is the total space
of a (real or complex) line bundle over the projective
model~$Y_{\AA,\GG}^{\P}$ which is isomorphic to the divisor component
in $Y_{\AA,\GG}$ indexed with $\bigcap \AA$.

The most prominent example of a projective arrangement model is the
minimal wonderful model for the complex braid arrangement,
$Y_{\AA_{n-2}^{\C}, \GG_{\rm min}}$. It is isomorphic to the
Deligne-Knudson-Mumford compactification $\overline{M_{0,n}}$ of the
moduli space of $n$-punctured complex projective
lines~\cite[4.3]{DP1}.

\medskip
Here is the first geometric interpretation of $D(\LL,\GG)$ in the case
of $\LL$ being the intersection lattice of a complex hyperplane arrangement.

\begin{thm}
\label{thm_geomint1}
{\rm (\cite{DP2, FY})} 
Let $\LL\,{=}\,\LL(\AA)$ be the intersection lattice of an essential 
arrangement  of complex hyperplanes $\AA$ and $\GG$ a building set in $\LL$
which contains the total intersection of $\AA$. Then, $D(\LL,\GG)$ is 
isomorphic to the integral cohomology algebra of the projective  
arrangement model~$Y_{\AA,\GG}^{\P}$:
\[
      D(\LL,\GG) \, \, \cong \, \, H^*(Y_{\AA,\GG}^{\P}, \Z)\, .
\] 
\end{thm}

\begin{explrm}
\label{ex_cohA2A3}
{\bf (Cohomology of braid arrangement models)} \newline
The projective arrangement model $Y_{\AA_2}^{\P}$ is homeomorphic to
the exceptional divisor in $Y_{\AA_2}\,{=}\,{\rm Bl}_{\{0\}}\C^2$, 
hence to $\C\P^1$. Its
cohomology is free of rank $1$ in degrees $0$ and $2$ and zero
otherwise.  Compare with $D(\Pi_3,\GG_{{\rm max}})$ in Example~\ref{ex_DPn}.

The projective arrangement model $Y_{\AA_3,\GG_{{\rm min}}}^{\P}$ is 
homeomorphic to $\overline{M_{0,5}}$, whose cohomology is known to be free of 
rank $1$ in 
degrees $0$ and $4$, free of rank $5$ in degree $2$, and zero otherwise. 
At least the coincidence of ranks is easy to verify in comparison with 
$D(\Pi_4,\GG_{{\rm min}})$ in Example~\ref{ex_DPn}.

Theorem~\ref{thm_geomint1} in fact gives an elegant presentation for 
the integral cohomology of 
$\overline{M_{0,n}}\,{\cong}\,Y_{\AA_{n-2},\GG_{{\rm min}}}^{\P}$ 
in terms of generators and relations:
\begin{eqnarray*}
  H^*(\overline{M_{0,n}}) &  \cong  &   D(\Pi_{n-1},\GG_{{\rm min}}) \\
                   &  \cong  &
       \Z\,[\,\{x_S\}_{S\subseteq [n-1], |S|\geq 2}\,] \Big/
        \left\langle \begin{array}{ll}
        x_S \, x_T & \mbox{for }\, 
                          S\cap T \neq \emptyset, \\
                & \mbox{and }\, 
                          S\not\subseteq T, T\not\subseteq S\, ,    \\[0.2cm]
        \sum_{\{i,j\}\subseteq S} \, x_S & 
                                          \mbox{for }\, 1\leq i<j\leq n-1\,
       \end{array}
       \right\rangle .     
\end{eqnarray*}

A lot of effort has been spent on describing the cohomology of 
$\overline{M_{0,n}}$ (cf~\cite{Ke}), none of the presentations comes 
close to the simplicity of the one stated above.

A nice expression for the  Hilbert function of  $H^*(\overline{M_{0,n}})$ 
has been derived by Yuzvinsky in \cite{Y1} as a consequence of his monomial 
linear basis for minimal projective arrangement models.
\end{explrm}

\medskip

To propose a more general geometric interpretation for $D(\LL,\GG)$, 
we start by describing a polyhedral fan $\Sigma(\LL,\GG)$ for 
any atomic lattice $\LL$ and any combinatorial building set $\GG$ in $\LL$.

\begin{dfrm} 
\label{df_Sigma}
{\bf (A simplicial fan realizing $\NN(\LL,\GG)$)} \newline
Let $\LL$ be an atomic lattice with set of atoms $\at=\{A_1,\ldots,A_n\}$,
$\GG$ a combinatorial building set in $\LL$.
For any $G\,{\in}\,\GG$ define the characteristic vector $v_G$ in~$\R^n$
by
\[
   (v_G)_i\, \, :=\, \, \left\{
\begin{array} {ll}
1 & \mbox{ if }\, G\geq A_i\, , \\
0 & \mbox{ otherwise\, , }
\end{array}
\right.   \qquad \mbox{ for }\, i=1,\ldots, n\, .
\] 
The simplicial fan $\Sigma(\LL,\GG)$ in $\R^n$ is the collection of cones
\[
       V_{\Ss}\,:=\,\mbox{cone} \{v_G\,|\, G\in \Ss\}
\]
for $\Ss$ nested in $\GG$.
\end{dfrm} 

By construction, $\Sigma(\LL,\GG)$ is a rational, simplicial fan that 
realizes the nested set complex $\NN(\LL,\GG)$. The fan gives rise to a
(non-compact) smooth toric variety $X_{\Sigma(\LL,\GG)}$~\cite[Prop.\ 2]{FY}.

\begin{explrm}
\label{ex_torvarPi3}
{\bf (The fan $\Sigma(\Pi_3,\GG_{{\rm max}})$ and its toric variety)}\newline
We depict $\Sigma(\Pi_3,\GG_{{\rm max}})$ in Figure~\ref{fig_fanPi3}.
The associated toric variety is the blowup of $\C^3$ in~$\{0\}$ with the proper 
transforms of coordinate axes removed. 

\begin{figure}[ht]
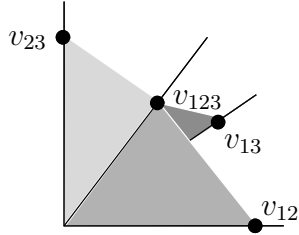

  \begin{picture}(0,0)%
    \includegraphics{fanPi3.pstex}%
  \end{picture}%
  \input{fanPi3.pstex_t}%
  
\caption{The simplicial fan $\Sigma(\Pi_3,\GG_{{\rm max}})$.} 
\label{fig_fanPi3}
\end{figure}
\end{explrm}

The algebra $D(\LL,\GG)$ here gains another geometric
meaning, this time for {\em any\/} atomic lattice~$\LL$. The abstract
algebraic detour of considering $D(\LL,\GG)$ in this general 
setting is rewarded by a somewhat unexpected return to geometry:

\begin{thm}
\label{thm_chow} {\rm \cite[Thm.\ 3]{FY}}
For an atomic lattice $\LL$ and a combinatorial building set $\GG$ 
in~$\LL$, $\DD(\LL,\GG)$ is isomorphic to the Chow ring of the toric 
variety $X_{\Sigma(\LL,\GG)}$,
\[
   D(\LL,\GG) \, \, \cong \, \, {\rm Ch}^*(X_{\Sigma(\LL,\GG)})\, .   
\]  
\end{thm}


\section{Adding arrangement models to the geometer's tool-box} 
\label{sect_gractns}

Let a diffeomorphic action of a finite group $\Gamma$ on a
smooth manifold $M$ be given. The goal is to modify the manifold by
blowups so as to have the group act on the resolution $\widetilde M$ 
with abelian stabilizers -- the quotient $\widetilde M/\Gamma$  
then has much more manageable singularities than the original quotient.
Such modifications for the sake of simplifying quotients have been 
of crucial importance at various places. One instance is Batyrev's work
on stringy Euler numbers~\cite{Ba1}, which in particular implies a 
conjecture of Reid~\cite{R}, and constitutes substantial progress 
towards higher dimensional MacKay correspondence.

There are two observations that point to wonderful arrangement models
as a possible tool in this context. First, the model construction is
equivariant if the initial setting carries a group action: if a
finite group $\Gamma$ acts on a real or complex 
vector space $V$, and the arrangement
$\AA$ is $\Gamma$-invariant, then the arrangement model $Y_{\AA,\GG}$
carries a natural $\Gamma$-action for any $\Gamma$-invariant building
set $\GG\subseteq \LL(\AA)$. Second, the model construction is not
bound to arrangements. In fact, locally finite stratifications of
manifolds which are local subspace arrangements, i.e., locally diffeomorphic 
to arrangements of linear subspaces, can be treated in a fully
analogous way. In the complex case, the construction has been pushed to
so-called {\em conical stratifications\/} by MacPherson and
Procesi~\cite{MP} with a real analogue developed by Gaiffi
in~\cite{G2}.
 
The significance of De~Concini-Procesi model constructions for
abelianizing group actions on complex varieties has been recognized by
Borisov and Gunnells~\cite{BG}, following work of
Batyrev~\cite{Ba1,Ba2}. Here we focus on the real setting.


\subsection{Learning from examples: permutation actions in low dimension}
\label{ssect_explactns}

Let us consider the action of the symmetric group $\mfS_n$ on real 
$n$-dimensional space by permuting coordinates:   
\[
    \sigma\,(x_1,\ldots, x_n)\, \, = \, \, 
    (x_{\sigma(1)},\ldots, x_{\sigma(n)}) \qquad 
    \mbox{for }\, \sigma \in \mfS_n, \, \,  x\in \R^n\, .
\]
Needless to say, we find a wealth of non-abelian stabilizers: For a
point $x\,{\in}\,\R^n$ that induces the set
partition $\pi=(B_1|\ldots|B_t)$ of $\{1,\ldots,n\}$ by pairwise 
coinciding coordinate entries, the stabilizer
of $x$ with respect to the permutation action is the Young subgroup
$\mfS_{\pi}\,{=}\,\mfS_{B_1} \,{\times}\, \ldots
\,{\times}\,\mfS_{B_t}$ of $\mfS_n$, where $\mfS_{B_i}$ denotes the
symmetric subgroup of $\mfS_n$ permuting the coordinates in $B_i$ for 
$i=1,\ldots,t$.

The locus of non-trivial stabilizers for the permutation action of
$\mfS_n$, in fact, is a familiar object: it is the rank~$n{-}1$ braid
arrangement $\AA_{n-1}$.  A natural idea that occurs when trying to 
abelianize a group action by blowups is to resolve the locus of {\em
non-abelian stabilizers\/} in a systematic way. Let us look at some
low dimensional examples.

\begin{explrm}
\label{ex_S3A2}
{\bf (The permutation action of $\mfS_3$)}\newline
We consider $\mfS_3$ acting on real $2$-space $V\,{\cong}\,\R^3/\Delta$.
The locus of non-trivial stabilizers consists of the $3$ hyperplanes
in $\AA_2$: for $x\,{\in}\,H_{ij}\,{\setminus}\,\{0\}$,
{\rm stab}$\,x\,{=}\,\langle (ij)\rangle\,{\cong}\, \Z_2$; in fact, 
$0$ is the only point having a non-abelian stabilizer, namely it is 
fixed by all of $\mfS_3$. 

Blowing up $\{0\}$ in $V$ according to the general idea outlined above, 
we recognize the maximal wonderful model for $\AA_2$ that we discussed 
in Example~\ref{ex_YA2}.

\begin{figure}[ht]
  \begin{picture}(0,0)%
    \includegraphics{S3A2.pstex}%
  \end{picture}%
  \input{S3A2.pstex_t}%
  
\caption{$\mfS_3$ acting on~$Y_{\AA_2}$.} \label{fig_S3A2}
\end{figure}

By construction, $\mfS_3$ acts coordinate-wise on $Y_{\AA_2}$.  
For points on proper transforms of hyperplanes $(y,H_{ij})\,{\in}\,D_{ij}$,
$1\,{\leq}\,i\,{<}\,j\,{\leq}\,3$, stabilizers are of order two:
stab$\,(y,H_{ij})\,{=}\,\langle(ij)\rangle\,{\cong}\, \Z_2$.
Otherwise, stabilizers are trivial, unless we are looking at one of
the three points $\psi_{ij}$ marked in Figure~\ref{fig_S3A2}. E.g.,
for $\psi_{12}{=}(0, \langle(1,-1,0)\rangle)$,
stab$\,\psi_{12}\,{=}\,\langle (12)\rangle\,{\cong}\,\Z_2$.  Although
the transposition $(12)$ does not fix the line
$\langle(1,-1,0)\rangle)$ point-wise, it fixes $\psi_{12}$ as a point
in $Y_{\AA_2}$!  We see that transpositions $(ij)\,{\in}\,\mfS_3$ act
on the open M\"obius band $Y_{\AA_2}$ by central symmetries in
$\psi_{ij}$.

Observe that the nested set stratification is not fine enough to
distinguish stabilizers: as the points $\psi_{ij}$ show, stabilizers 
are not isomorphic on open strata. 
\end{explrm}

\begin{explrm}
\label{ex_S4A3}
{\bf (The permutation action of $\mfS_4$)}\newline
Let us now consider $\mfS_4$ acting on real $3$-space 
$V\,{\cong}\,\R^4/\Delta$.
The locus of non-abelian stabilizers consists of the triple 
intersections of hyperplanes in $\AA_3$, i.e., the subspaces contained 
in the minimal building set $\GG_{{\rm min}}$ in $\LL(\AA_3){=}\Pi_4$.
Our general strategy suggests to look at the arrangement model 
$Y_{\AA,\GG_{{\rm min}}}$.

\begin{figure}[ht]
  \begin{picture}(0,0)%
    \includegraphics{S4A3.pstex}%
  \end{picture}%
  \input{S4A3.pstex_t}%
  
\caption{$\mfS_4$ acting on~${\rm Bl}_{\{0\}}V$,
$V\,{=}\,\R^4{/}\Delta$.} \label{fig_S4A3}
\end{figure}

We consider a situation familiar to us from Example~\ref{ex_A3}. In
Figure~\ref{fig_S4A3}, we illustrate the situation after the first
blowup step in the construction of $Y_{\AA,\GG_{{\rm min}}}$, i.e.,
the exceptional divisor after blowing up $\{0\}$ in $V$ with the
stratification induced by the hyperplanes of $\AA_3$. To complete the
construction of $Y_{\AA,\GG_{{\rm min}}}$, another $4$ blowups in the
triple intersections of hyperplanes are necessary, the result of which
we illustrate locally for the triple intersection corresponding
to~$134$.
Triple intersections of hyperplanes in ${\rm Bl}_{\{0\}}V$ have
stabilizers isomorphic to $\mfS_3$ -- the further blowups in
triple intersections are indeed necessary to obtain an abelianization
of the permutation action.

Again, we observe that the nested  set 
stratification on $Y_{\AA,\GG_{{\rm min}}}$
does not distinguish stabilizers: we indicate subdivisions of nested set
strata resulting from non-isomorphic stabilizers by dotted lines, 
respectively unfilled points in Figure~\ref{fig_S4A3}.

Let us look at stabilizers of points on the model $Y_{\AA,\GG_{{\rm min}}}$:
We find points with stabilizers isomorphic to $\Z_2$ -- any 
generic point on a divisor $D_{ij}$ will be such. We also find points 
with stabilizers isomorphic to $\Z_2\,{\times}\,\Z_2$, e.g., the point 
$\omega$ on 
$D_{1234}$ corresponding to the line $\langle(1,-1,0,0)\rangle$. 

But, on $Y_{\AA,\GG_{{\rm min}}}$ we also find points with {\em
non-abelian\/} stabilizers! For example, the intersection of $D_{14}$ and 
$D_{23}$ on $D_{1234}$ corresponding to the line 
$\langle(1,-1,-1,1)\rangle$ is stabilized by 
both $(14)$ and $(12)(34)$ in $\mfS_4$, which do not commute. In fact,
the stabilizer is isomorphic to $\Z_2\,{\wr}\,\Z_2$. 

This observation shows that blowing up the locus of non-abelian
stabilizers is not enough to abelianize the action! Further blowups in
double intersections of hyperplanes are necessary, which suggests,
contrary to our first assumption, the {\em maximal\/} arrangement
model $Y_{\AA_3}$ as an abelianization of the permutation action.

Some last remarks on this example: observe that stabilizers of points 
on~$Y_{\AA_3}$ all are elementary abelian $2$-groups. We will later see that
the strategy of resolving finite group actions on real vector spaces and 
even manifolds by constructing a suitable maximal De~Concini-Procesi 
model does not only abelianize the action, but yields stabilizers
isomorphic to elementary abelian $2$-groups.

Also, it seems we cannot do any better than that within the framework
of blowups, i.e., we neither can get rid of non-trivial stabilizers,
nor can we reduce the rank of non-trivial stabilizers any further. The
divisors $D_{ij}$ are stabilized by transpositions $(ij)$ which supports
our first claim.  For the second claim, consider the point 
$\omega\,{=}\,(0,\ell_1)$ 
in $Y_{\AA_3}$ with $\ell_1\,{=}\,\langle(1,-1,0,0)\rangle$ 
(here we use the encoding 
of points on arrangement models proposed in~(\ref{eq_encoding})).  
We have seen above that
stab$\,\omega\,{\cong}\,\Z_2\,{\times}\,\Z_2$, in fact
stab$\,\omega\,{=}\,\langle(12)\rangle\,{\times}\,\langle(34)\rangle$.
Blowing up $Y_{\AA_3}$ in~$\omega$ means to again glue in an open M\"obius
band. Points on the new exceptional divisor $D_{\omega}\,{\cong}\,\R\P^1$ will be
parameterized by tupels $(0,\ell_1,\ell_2)$, where 
$\ell_2$ is a line orthogonal to $\ell_1$ in $V$. A
generic point on this stratum will be stabilized only by the 
transposition~$(12)$,
specific points however, e.g.,
$(0,\ell_1,\langle(0,0,1,-1)\rangle)$ will still be
stabilized by all of stab$\,\omega\,{\cong}\,\Z_2\,{\times}\,\Z_2$.
\end{explrm}


\subsection{Abelianizing a finite linear action}
\label{ssect_linactns}

Following the basic idea of proposing De~Concini-Procesi arrangement
models as abelianizations of finite group actions and drawing from our
experiences with the permutation action on low-dimensional real space 
in Section~\ref{ssect_explactns} we here treat the case of finite linear 
actions.

Let a finite group $\Gamma$ act linearly and effectively on
real $n$-space~$\R^n$.  Without loss of generality, we can assume
that the action is orthogonal~\cite[2.3, Thm.~1]{V}; we fix the
appropriate scalar product throughout.

Our strategy is to construct an arrangement of subspaces $\AA(\Gamma)$
in real $n$-space, and to propose the maximal wonderful model 
$Y_{\AA(\Gamma)}$ as an abelianization of the given action.

\begin{cstrrm} 
\label{cstr_arrgt}
{\bf (The arrangement $\AA(\Gamma)$)} \newline
For any subgroup $H$ in $\Gamma$, define a linear subspace
\begin{equation} \label{eq_Lspaces}
    L(H)\, \, := \, \, 
{\rm span}\{\,\ell\,|\, 
\ell \mbox{ line in }\, \R^n \mbox{ with }\, 
            H\circ \ell = \ell\,\}\, , 
\end{equation}
the linear span of all lines in $V$ that are invariant under the 
action of~$H$.

Denote by $\AA(\Gamma)\,=\,\AA(\Gamma \circlearrowright \R^n)$ 
the arrangement of proper subspaces in $\R^n$
that are of the form $L(H)$ for some subgroup $H$ in $\Gamma$.
\end{cstrrm}

Observe that the arrangement $\AA(\Gamma)$ never contains any
hyperplane: if $L(H)$ were a hyperplane for some subgroup $H$ in
$\Gamma$, then also its orthogonal line $\ell$ would be
invariant under the action of $H$. By definition of $L(H)$, however,
$\ell$ would then be contained in $L(H)$ which in turn would be the full
ambient space.

\begin{thm}\label{thm_linactn} {\rm \cite[Thm.\ 3.1]{FK3}}
For any effective linear action of a finite group~$\Gamma$ on
$n$-dimensional real space, the maximal wonderful arrangement model
$Y_{\AA(\Gamma)}$ abelianizes the action. Moreover, stabilizers of
points on the arrangement model are isomorphic to 
elementary abelian $2$-groups.
\end{thm}

The first example coming to mind is the permutation action of $\mfS_n$
on real $n$-space. We find that $\AA(\mfS_n)$ is the rank~$2$
truncation of the braid arrangement, $\AA_{n-1}^{{\rm rk}\geq 2}$, i.e.,
the arrangement consisting of subspaces in $\AA_{n-1}$ of 
codimension~${\geq}\,2$.
For details, see~\cite[Sect. 4.2]{FK3}. In earlier work~\cite{FK2}, we
had already proposed the maximal arrangement model of the braid
arrangement as an abelianization of the permutation action. 
We proved that stabilizers on $Y_{\AA_{n-1}}$ are isomorphic to 
elementary abelian $2$-groups by providing 
explicit descriptions of stabilizers based on an algebraic-combinatorial 
set-up for studying these groups.

\subsection{Abelianizing finite diffeomorphic actions on manifolds}
\label{ssect_actnsmfds}

Let us now look at a generalization of the abelianization presented in
Section~\ref{ssect_linactns}. Assume that $\Gamma$ is a finite group
that acts diffeomorphically and effectively on a smooth real
manifold~$M$.  We first observe that such an action induces a linear
action of the stabilizer stab$\, x$ on the tangent space $T_xM$ at any
point $x$ in $M$. Hence, locally we are back to the setting that we 
discussed before: For any subgroup $H$ in stab$\, x$, we can define a 
linear subspace $L(x,H):=L(H)$ of the tangent space $T_xM$ as 
in~(\ref{eq_Lspaces}), and we can combine the non-trivial subspaces to
form an arrangement 
$\AA_{x}\,{:=}\,\AA({\rm stab}\, \circlearrowright T_xM)$ in $T_xM$.

Combined with the information that a model construction in the spirit
of De~Concini-Procesi arrangement models exists also for local
subspace arrangements, we need to stratify the manifold so as to
locally reproduce the arrangement $\AA_{x}$ in any tangent
space $T_xM$. Here is how to do that:

\begin{cstrrm}
\label{cstr_stratt}
{\bf (The stratification $\mathfrak L$)} \newline
For any $x\,{\in}\, M$, and any subgroup $H$ in stab$\,x$, define 
a normal (!) subgroup $F(x,H)$ in $H$ by
\[
       F(x,H)\, \, = \, \, \{h\in H\,|\, h \circ y=y \mbox{ for any }
\, y\in L(x,H) \}\, ;
\]
$F(x,H)$ is the subgroup of elements in $H$ that fix all of $L(x,H)$ 
point-wise. Define $\mfL (x,H)$ to be the connected component of the 
fixed point set of $F(x,H)$ in $M$ that contains~$x$. Now combine these 
submanifolds so as to form a locally finite stratification
\[
      \mfL\, \, = \, \,  (\mfL(x,H))_{\,
                              x\in M, \, 
                              H\leq {\rm stab}\, x} \, .
\] 
\end{cstrrm}

Observe that, as we tacitly did for stratifications induced by arrangements 
or by irreducible components of divisors, we only specify strata of 
proper codimension. 

The stratification $\LL$ locally coincides with the tangent space stratifications 
coming from our linear setting. Technically speaking: for any $x\,{\in}\,M$,
there exists an open neighborhood $U$ of $x$ in $M$, and a 
stab$\,x$-equivariant diffeomorphism $\Phi_x:\, U \rightarrow 
T_xM$ such that 
\begin{equation} \label{eq_locsubsparrgt}
\Phi_x(\mfL(x,H))\,{=}\,L(x,H)\,
\end{equation} 
for any subgroup $H$ in stab$\,x$.  In particular,
(\ref{eq_locsubsparrgt}) shows that the stratification $\mfL$ of $M$
is a local subspace arrangement.

\begin{thm}\label{thm_actnmfd} {\rm \cite[Thm.\ 3.4]{FK3}}
Let a finite group $\Gamma$ act diffeomorphically and effectively on
a smooth real manifold $M$. Then the wonderful model $Y_{\mfL}$ induced 
by the locally finite stratification $\mfL$ of $M$ abelianizes the action. 
Moreover, stabilizers of points on the model $Y_{\mfL}$  are isomorphic 
to elementary abelian $2$-groups.
\end{thm}

\begin{explrm}
\label{ex_S3RP2}
{\bf (Abelianizing the permutation action on $\R\P^2$)}\newline
Let us look at a small non-linear example: the permutation action 
of $\mfS_3$ on the real projective plane induced by $\mfS_3$ permuting 
coordinates in $\R^3$. 

We picture $\R\P^2$ by its upper hemisphere model in 
Figure~\ref{fig_S3RP2}, where we agree to place the projectivization of 
$\Delta^{\perp}$ on the equator. The locus of non-trivial stabilizers of the 
$\mfS_3$ permutation action consists of the projectivizations of 
hyperplanes $H_{ij}{:}\, x_i\,{=}\,x_j$, $1\,{\leq}\,i\,{<}\,j\,{\leq}\,3$,
and three additional points $\Psi_{ij}$ on $\P\Delta^{\perp}$ indicated 
in Figure~\ref{fig_S3RP2}. The $\mfS_3$ action can be visualized by 
observing that transpositions $(ij)\,{\in}\, \mfS_3$ 
act as reflections in the lines $\P H_{ij}$, respectively.

\begin{figure}[ht]
  \begin{picture}(0,0)%
    \includegraphics{S3RP2.pstex}%
  \end{picture}%
  \input{S3RP2.pstex_t}%
  
\caption{$\mfS_3$ acting on~$\R\P^2$: the stabilizer stratification.} 
\label{fig_S3RP2}
\end{figure}

We find that the arrangements $\AA_{\ell}$ in the tangent spaces
$T_{\ell} \R\P^2$ are empty, unless
$\ell\,{=}\,[1{:}1{:}1]$. Hence, (\ref{eq_locsubsparrgt}) allows us 
to conclude that the $\mfL$-stratification of $\R\P^2$ consists 
of a single point, $[1{:}1{:}1]$. Observe that the $\mfS_3$-action on 
$T_{[1:1:1]} \R\P^2$ coincides with the permutation action of 
$\mfS_3$ on $\R^3/\Delta$. 

The wonderful model $Y_{\mfL}$ hence is a Klein bottle, the result 
of blowing up $\R\P^2$ in $[1{:}1{:}1]$, i.e., glueing a M\"obius 
band into the punctured projective plane. 

Observe that the $\mfL$-stratification is coarser than the codimension $2$
truncation of the stabilizer stratification: 
The isolated points $\Psi_{ij}$ on $\P\Delta^{\perp}$ have 
non-trivial stabilizers, but do not occur as strata in 
the $\mfL$-stratification. 
\end{explrm}


\end{document}